\documentclass{article}

\usepackage{arxiv}

\usepackage[utf8]{inputenc} % allow utf-8 input
\usepackage[T1]{fontenc}    % use 8-bit T1 fonts
\usepackage{hyperref}       % hyperlinks
\usepackage{url}            % simple URL typesetting
\usepackage{booktabs}       % professional-quality tables
\usepackage{amsfonts}       % blackboard math symbols
\usepackage{nicefrac}       % compact symbols for 1/2, etc.
\usepackage{microtype}      % microtypography
\usepackage{lipsum}
\usepackage{graphicx}
\usepackage{amsmath}
\usepackage{amsfonts}
\usepackage{amssymb}
\usepackage{subfigure}
\usepackage{multirow}
\usepackage{float}
\usepackage{adjustbox}
\usepackage{booktabs}
\usepackage{lscape}
\usepackage{enumitem, array}

\title{DATA-DRIVEN MODELING OF SHIP MANEUVERS IN WAVES VIA DYNAMIC MODE DECOMPOSITION}

\author{
 Matteo Diez$^{1,\star}$, Andrea Serani$^{1}$, Emilio F. Campana$^2$, and Frederick Stern$^3$\\
 \\
  $^1$CNR-INM, National Research Council-Institute of Marine Engineering, Rome, Italy\\
  $^2$ CNR-DIITET, National Research Council\\
  Engineering, ICT and Technologies for Energy and Transportation Department, Rome, Italy\\
  $^3$IIHR-Hydroscience \& Engineering, The University of Iowa, Iowa City, USA\\
  $^\star$Corresponding author: \texttt{matteo.diez@cnr.it} \\
}

\begin{document}
\maketitle

\begin{abstract}
A data-driven and equation-free approach is proposed and discussed to model ships maneuvers in waves, based on the dynamic mode decomposition (DMD). 
DMD is a dimensionality-reduction/reduced-order modeling method, which provides a linear finite-dimensional representation of a possibly nonlinear system dynamics by means of a set of modes with associated oscillation frequencies and decay/growth rates. DMD also allows for short-term future estimates of the system's state, which can be used for real-time prediction and control.
Here, the objective of the DMD is the analysis and forecast of the trajectories/motions/forces of ships operating in waves, offering a complementary efficient method to equation-based system identification approaches. 
Results are presented for the course keeping of a free-running naval destroyer (5415M) in irregular stern-quartering waves and for the free-running KRISO Container Ship (KCS) performing a turning circle in regular waves. Results are overall promising and show how DMD is able to identify the most important modes and forecast the system's state with reasonable accuracy upto two wave encounter periods.
\end{abstract}

% keywords can be removed

\section{Introduction}
To ensure the safety of structures, payload, and crew in adverse weather conditions, ships must have good seakeeping, maneuverability, and structural performance. In this regard, commercial and military ships must meet International Maritime Organization (IMO) Guidelines and NATO Standardization Agreements (STANAG). The prediction capability of ship performance in waves, along with the understanding of the physics involved, is of utmost importance. 
Recent NATO Science and Technology Organization (STO) Applied Vehicle Technology (AVT) task groups, such as AVT-280 ``Evaluation of Prediction Methods for Ship Performance in Heavy Weather'' (2017-2019) and AVT-348 ``Assessment of Experiments and Prediction Methods for Naval Ships Maneuvering in Waves'' (2021-2023) focus on prediction methods for ship seakeeping and maneuvering in waves, respectively.
Recent computational and experimental fluid dynamics studies have demonstrated the maturity of computational tools for the prediction of ship performance in waves, including their assessment in extreme sea conditions \cite{van2020prediction,serani2021-OE}. The computational cost associated with the analysis is generally very high, especially if statistical convergence of relevant estimators is sought after and complex hydro-structural problems are investigated via high-fidelity solvers \cite{diez2020fsi}. In this context, machine learning and reduced-order models can help reducing the computational cost providing fast estimates, given that they are properly trained and/or calibrated. In addition, the ability of these methods to learn from data makes them suitable for use in digital twins platforms. Finally, reduced-order models are generally easier to interpret than machine learning approaches and could help shedding light on the physics involved.

%Computations and experiments usually produce a large amount of data, whose investigation could shed light on the problem physics. 
%The use of reduced-order models generally reduces the computational cost and could help shedding light on the physical phenomena involved.

Here, we propose and discuss a data-driven and equation-free approach to the reduced-order modeling of ships maneuvering in waves, based on the dynamic mode decomposition (DMD). 
DMD is a dimensionality-reduction/reduced-order modeling method, which provides a set of modes with associated oscillation frequencies and decay/growth rates \cite{schmid2010dynamic}. For linear systems, these modes/frequencies correspond to the linear normal modes/frequencies of the system. More generally, DMD modes/frequencies approximate eigenmodes and eigenvalues of the infinite-dimensional linear Koopman operator, providing a linear finite-dimensional representation of the (possibly nonlinear) system dynamics \cite{kutz2016multiresolution}.
Its growing success, especially in the fluid dynamics community (e.g. \cite{rowley2009spectral,magionesi2018modal,timur2020}), is due to its equation-free and data-driven nature. The method is capable of providing accurate assessments of the spatio-temporal coherent structures in complex flows and systems, also allowing short-term future estimates of the system's state, which can be used for real-time prediction and control \cite{kutz2016dynamic}.

In the present work, the objective of the DMD is the %data-driven and equation-free 
analysis and forecast of the finite-dimensional set of trajectory/motion/force time histories of ships operating in waves, offering a complementary efficient method to equation-based system identification approaches, e.g.  \cite{araki2012estimating,araki2019improved}. The efficiency of the method in this context stems from the finite dimensionality of the set of relevant state variables together with the simplicity of operations required to model the system dynamics (as opposed to more data/resource-consuming machine learning approaches). This offers opportunity for integration into digital twin platforms for the data-driven modeling and prediction of ships in waves. 

Results are presented for the course keeping of a free-running naval destroyer (5415M) in irregular stern-quartering waves at target $\mathrm{Fr}= 0.33$ and sea state 7, using URANS (unsteady Reynolds-averaged Navier-Stokes) computations from \cite{serani2021-OE}.
Results are also presented for the free-running KRISO Container Ship (KCS), using experimental data from the University of Iowa IIHR wave basin, focusing on starboard turning circle with rudder angle of 35\textit{deg}, target Froude number $\mathrm{Fr}=0.157$ in regular waves with $\lambda/L=1$ (wavelength to ship-length ratio) and $H/\lambda=1/60$ (wave-height to wavelength ratio). The present research is conducted in collaboration with NATO AVT-348 ``Assessment of Experiments and Prediction Methods for Naval Ships Maneuvering in Waves'', and AVT-351 ``Enhanced Computational Performance and Stability \& Control Prediction for NATO Military Vehicles''.

\section{Dynamic Mode Decomposition}
The DMD formulation and nomenclature is taken from \cite{kutz2016dynamic}. Specifically, consider a dynamical system described as
\begin{equation}\label{eq:sysdyn}
\dfrac{\mathrm{d}\mathbf{x}}{\mathrm{d}t}=\mathbf{f}(\mathbf{x},t;\mu),
\end{equation}
where $\mathbf{x}(t)\in\mathbb{R}^n$ represents the system's state at time $t$, $\mu$ contains the parameters of the system, and $\mathbf{f}(\cdot)$ represents its dynamics.
% In general, a dynamical system can be described by ordinary differential equations (often nonlinear). 
The state $\mathbf{x}$ is generally large, with $n\gg 1$ and can represents, for instance, the discretization of partial differential equations at a number of discrete spatial points, or multi-channel/multi-variable time series.

%Generally, it is not possible to find a solution for a nonlinear dynamical system, so numerical methods are used to evolve to future states. 
Considering $\mathbf{f}(\mathbf{x},t;\mu)$ as unknown, the DMD works with an equation-free perspective. Thus, only the system measurements are used to approximate the system dynamics and forecast the future states. Equation \ref{eq:sysdyn} is approximated by the DMD as a locally linear dynamical system defined as
\begin{equation}\label{eq:dmdsys}
\dfrac{\mathrm{d}\mathbf{x}}{\mathrm{d}t}=\mathcal{A}\mathbf{x}
\end{equation}
with solution
\begin{equation}\label{eq:dmdrec}
\mathbf{x}(t)=\sum_{k=1}^n \boldsymbol{\phi}_k \, q_k(t)=
\sum_{k=1}^n \boldsymbol{\phi}_k \, b_k \,\exp(\omega_k t),
\end{equation}
where $\boldsymbol{\phi}_k$ and $\omega_k$ are, respectively, the eigenvectors and the eigenvalues of the matrix $\mathcal{A}$ and the coefficients $b_k$ are the coordinates of the initial condition $\mathbf{x}_0$ in the eigenvector basis, i.e. $\mathbf{b}=\boldsymbol{\Phi}^{-1}\mathbf{x}_0$.

%Discretizing the time dynamics in Eq. \ref{eq:sysdyn} by 
Sampling the system every $\Delta t$, the time-discrete state can be expressed as $\mathbf{x}_k=\mathbf{x}(k\Delta t)$ with $k=1,...,m$, representing from now on known system measurements.
An equivalent discrete-time representation of the system in Eq. \ref{eq:dmdsys} can be written as
\begin{equation}\label{eq:dmddsys}
\mathbf{x}_{k+1}=\mathbf{Ax}_k, \qquad \mathrm{with} 
\qquad \mathbf{A}=\exp(\mathcal{A}\Delta t)
\end{equation}

Arranging all the $m$ system measurements in the following two matrices
\begin{equation}
\mathbf{X}=
\begin{bmatrix}
| & | & & |\\
\mathbf{x}_1 & \mathbf{x}_2 & \dots & \mathbf{x}_{m-1}\\
| & | & & |\\
\end{bmatrix},
\qquad
\mathbf{X}'=
\begin{bmatrix}
| & | & & |\\
\mathbf{x}_2 & \mathbf{x}_3 & \dots & \mathbf{x}_{m}\\
| & | & & |\\
\end{bmatrix},
\end{equation}
the matrix $\mathbf{A}$ in Eq. \ref{eq:dmddsys} can be constructed using the following approximation
\begin{equation}\label{eq:approxA}
\mathbf{A}\approx\mathbf{X}'\mathbf{X}^{\dag},
\end{equation}
where $\mathbf{X}^{\dag}$ is the Moore-Penrose pseudoinverse of $\mathbf{X}$, which minimize $\|\mathbf{X}'-\mathbf{AX}\|_F$, where $\|\cdot\|_F$ is the Frobenius norm.

%The matrix $\mathbf{A}$ can be constructed from the system measurements using the following approximation
%%
%\begin{equation}\label{eq:approxA}
%\mathbf{A}\approx\mathbf{X}'\mathbf{X}^{\dag},
%\end{equation}
%%
%where 
%\begin{equation}
%\mathbf{X}=
%\begin{bmatrix}
%| & | & & |\\
%\mathbf{x}_1 & \mathbf{x}_2 & \dots & \mathbf{x}_{m-1}\\
%| & | & & |\\
%\end{bmatrix},
%\qquad
%\mathbf{X}'=
%\begin{bmatrix}
%| & | & & |\\
%\mathbf{x}_2 & \mathbf{x}_3 & \dots & \mathbf{x}_{m}\\
%| & | & & |\\
%\end{bmatrix},
%\end{equation}
%%
%and $\mathbf{X}^{\dag}$ is the Moore-Penrose pseudoinverse of $\mathbf{X}$, which minimize $\|\mathbf{X}'-\mathbf{AX}\|_F,$, where $\|\cdot\|_F$ is the Frobenius norm. 

The state-variable evolution in time can be approximated by the same modal expansion of Eq. \ref{eq:dmdrec}, where the $\boldsymbol\phi_k$ are now the eigenvectors of the approximated matrix $\bf A$ and $\omega_k=\text{ln}(\lambda_k)/\Delta t$ with $\lambda_k$ eigenvalues of the same matrix \cite{kutz2016dynamic}.

In general, the DMD can be viewed as a method to compute the eigenvalues and eingevectors (modes) of a finite-dimensional linear model that approximates the infinite-dimensional linear Koopman operator \cite{kutz2016dynamic}, also known as the composition operator. Here, the DMD is applied to inherently finite-dimensional data, i.e. ship trajectories/motions/forces in waves, similarly to DMD applications to power grid load data \cite{dylewsky2020dynamic,mohan2018data}, financial trading strategies \cite{mann2016dynamic}, sales data \cite{vasconcelos2019dynamic}, and neural recordings \cite{brunton2016extracting}.
Furthermore, due to the low dimensionality of data in the current context, Eq. \ref{eq:approxA} is computed directly, without the need of performing the singular value decomposition of $\bf X$ and projecting onto proper orthogonal decomposition modes \cite{kutz2016dynamic}.

\section{Test Cases and DMD Setup}
This section describes the test cases used for demonstration and the DMD setup. It is worth noting that, although the DMD is not a machine learning method in the strict sense, its data-driven nature allows for approaching DMD in a similar way to machine learning. 
Here the matrix $\bf A$ is constructed using observed (past) time histories, which are used as training set. The DMD is then used for the short-term prediction of trajectory/motion/force time histories, which are compared against true observed (future) data used as test set. As comparison metrics, the average normalized mean square error (NMSE) is used. The mean square of the modal coordinates $\langle q_k^2 \rangle$ is used as a metrics for modal participation.
\begin{figure}[!b]
\centering
\includegraphics[width=0.45\textwidth]{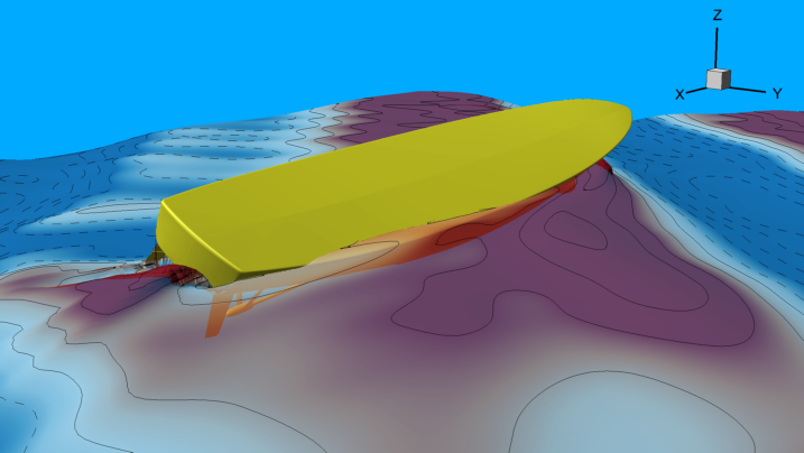}\hspace{0.5cm}
\includegraphics[width=0.45\textwidth]{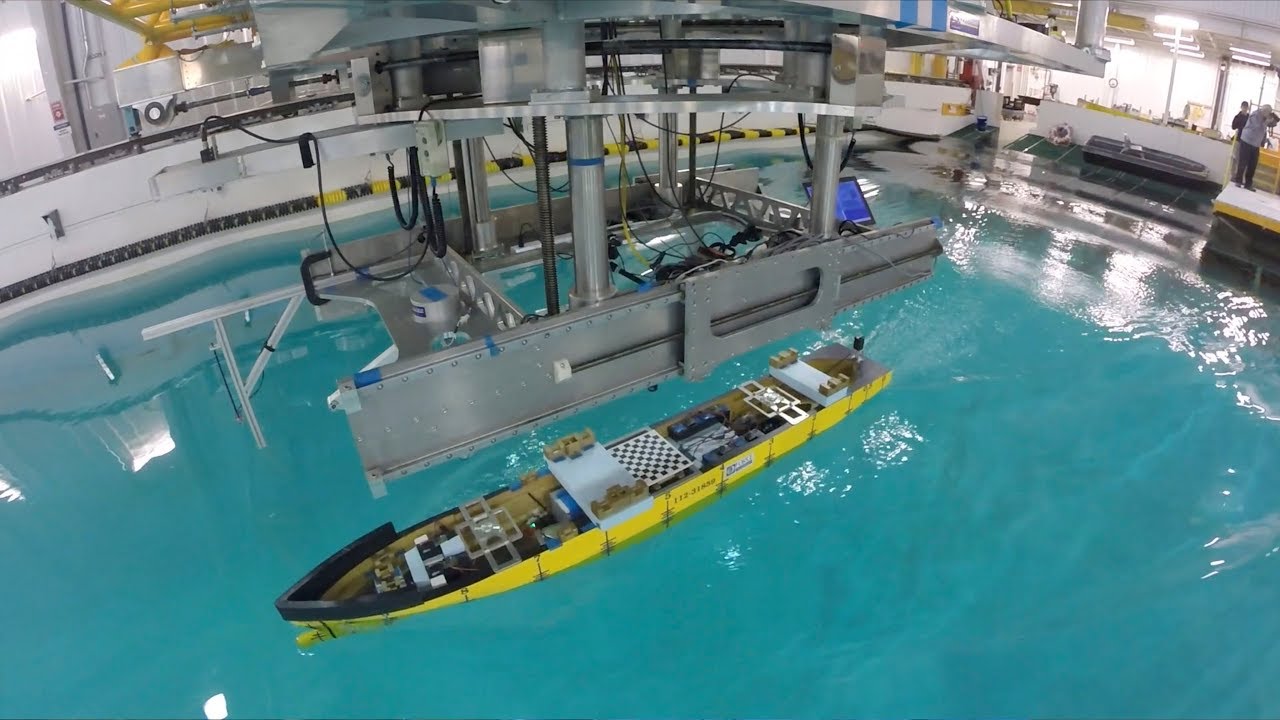}
\caption{CFD snapshot of the 5415M test case (left) and a photograh of the IIHR wave basin with the KCS model (right).}\label{fig:5415}
\end{figure}

\subsection{Course Keeping of the 5415M in Irregular Waves}
The hull form under investigation is the MARIN model 7967 which is equivalent to 5415M, used as test case for the NATO AVT-280 ``Evaluation of Prediction Methods for Ship Performance in Heavy Weather'' \cite{van2020prediction}. The model is self-propelled and kept on course by a proportional-derivative controller actuating the rudder angle. 

Course-keeping computations are based on the URANS code CFDShip-Iowa V4.5 \cite{huang2008-IJNMF}. CFD simulations are performed with propeller revolutions per minute fixed to the self-propulsion point of the model for the nominal speed, corresponding to $\mathrm{Fr} = 0.33$. The simulations are conducted in irregular long crested waves (following a JONSWAP spectrum), with nominal peak period $T_p = 9.2$\textit{s} and wave heading of 300\textit{deg} (see Fig. \ref{fig:5415}). The nominal significant wave height is equal to 7\textit{m}, corresponding to sea state 7 (high), according to the World Meteorological Organization definition. It may be noted that the simulation conditions are close to a resonance condition for the roll.
%
%
%The turbulence is computed by the isotropic Menter's blended $k-\epsilon/k-\omega$ (BKW) model with shear stress transport (SST) using no wall function. The location of the free surface is given by the ''zero'' value of the level-set function, positive in the water and negative in the air. 
The six degrees of freedom rigid body equations of motion are solved to calculate linear and angular motions of the ship. A simplified body-force model is used for the propeller, which prescribes axisymmetric body force with axial and tangential components. The total number of grid points is about 45M. Further details can be found in \cite{serani2021-OE} and \cite{van2020prediction} where also potential flow computations and experimental data are presented and discussed. 

The state variables used for DMD are the ship six degrees of freedom (surge, sway, heave, roll, pitch, yaw; ship motions in the carriage coordinate system, projected onto the ship axes) plus the rudder angle. In addition, first and second time derivatives of all variables are included in the data set, computed by a fourth-order finite difference scheme. The use of derivatives enables a better description of the system dynamics and defines a higher-dimensional space potentially more amenable to an accurate linear representation. All variables are standardized, i.e. translated and scaled to have zero mean and unit variance. The matrix $\bf A$ is built using about five encounter waves (corresponding to 1766 time steps). The prediction and test sets span the same length. Finally, the prediction is built using all modes/frequencies, leaving reduced-order studies to future research.

\subsection{Turning Circle of the KCS in Regular Waves}

The second test case is the starboard turning circle of the free-running KCS in regular waves with constant rudder angle of 35\textit{deg}. Data are taken from experiments conducted at the IIHR wave basin, which is shown in Fig. \ref{fig:5415} and whose characteristics are given in \cite{sanada2021experimental}. The model length is $L=2.7$\textit{m} and the nominal speed corresponds to $\mathrm{Fr}=0.157$. The propeller RPM are fixed and provide the nominal speed in calm water. The regular wave parameters are $\lambda/L=1$ and $H/\lambda=1/60$.

The state variables used for DMD are the $x$, $y$, and $z$ coordinates (Earth-coordinate system) of a reference point placed amidships, pitch and roll motions, turning rate, $u$, $v$, and $w$ component of the ship velocity (projected onto the ship axes), rudder angle, propeller thrust, and torque. As in the previous test case, first and second time derivatives are included in the data set and all variables are standardized. The matrix $\bf A$ is built using about four encounter waves (corresponding to 132 time steps), while the prediction and test sets span the same length. Predictions are build using all modes/frequencies.

\section{Results}
The following subsections describe the DMD results for the test cases, discussing the DMD analysis of the system dynamics (including complex modal frequencies, modal participation, and most energetic modes) and the future prediction of the system's state.

\begin{figure}[!b]
\centering
\includegraphics[width=1\textwidth]{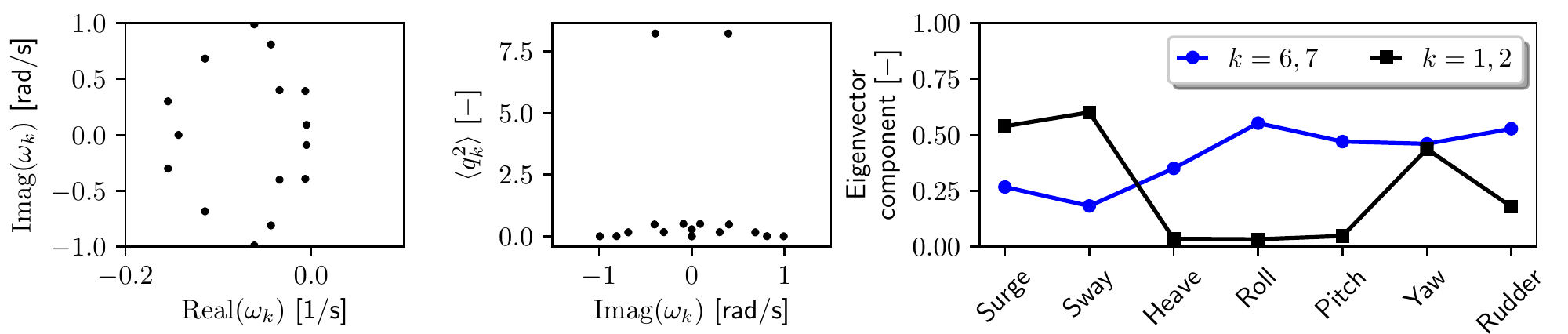}
\caption{DMD complex modal frequencies, modal participation, and two most energetic modes (from left to right, respectively) for the 5415M test case.}\label{fig:freq_5415}
\end{figure}
\begin{figure}[!t]
\centering
\includegraphics[width=0.975\textwidth]{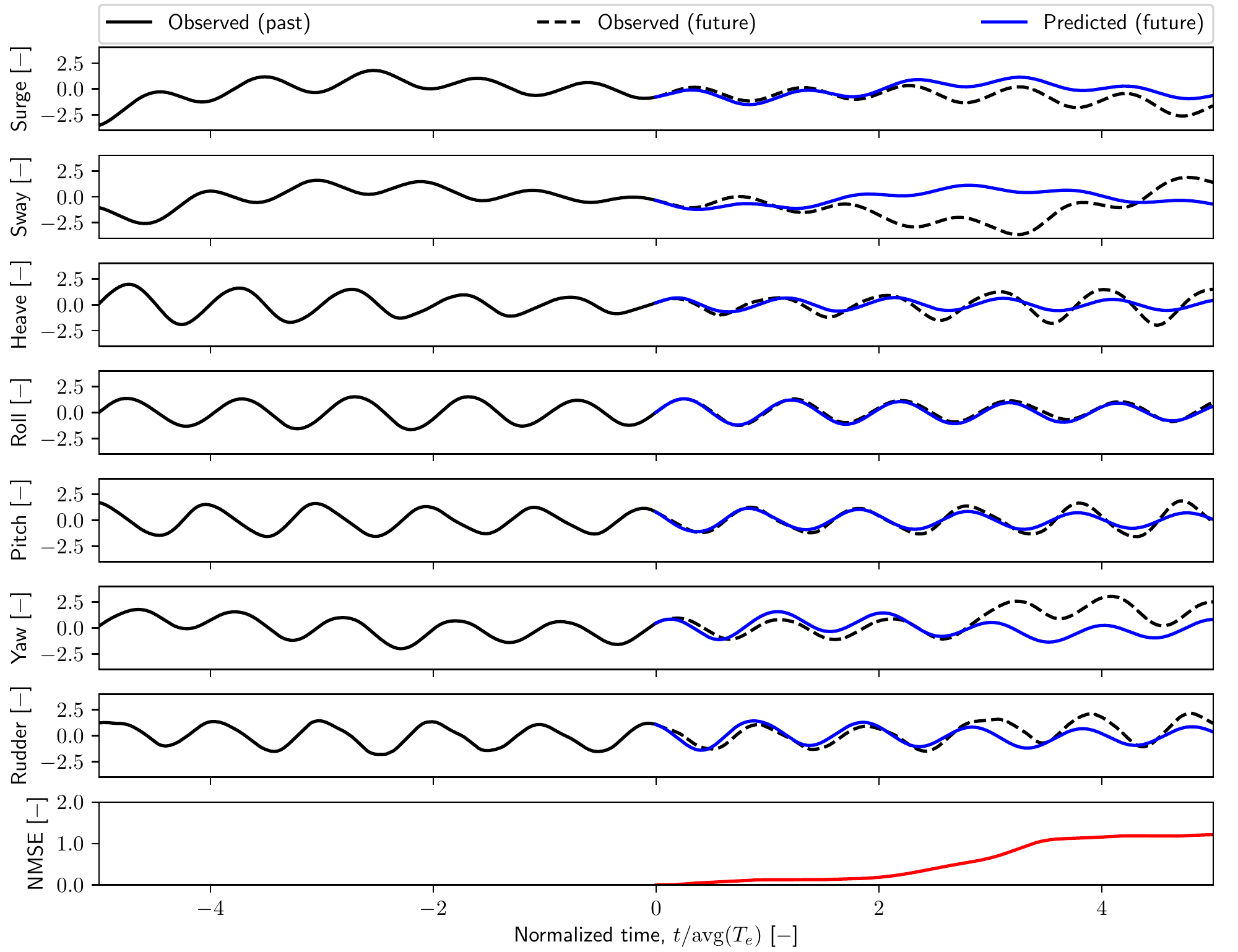}
\caption{DMD prediction and average error for the 5415M test case (standardized variables).}\label{fig:pred_5415}
\end{figure}

\subsection{5415M Course Keeping in Irregular Waves}
Figure \ref{fig:freq_5415} presents the DMD results for the system dynamics. Specifically, the left figure shows the complex modal frequencies provided by DMD. The center figure shows the modal participation as a function of the frequency imaginary part. Finally, the right figure presents the magnitude of components for the two most energetic couples of complex-conjugate modes, where for the sake of clarity and simplicity derivatives are not shown. It may be noted how the dynamic is largely participated by one couple only (see Fig. \ref{fig:freq_5415} center, namely represented by the modes $k=6,7$). This couple presents a frequency close to the roll resonance frequency (see, e.g., \cite{van2020prediction}) and mainly involves roll, pitch, yaw, rudder, and to a lesser extent heave (see blue line in Fig. \ref{fig:freq_5415} right). The second couple is significantly less energetic and has a smaller frequency (slower dynamics). It mainly involves planar motion variables, i.e. surge, sway, yaw, and to a lesser extent rudder (see black line in Fig. \ref{fig:freq_5415} right).

The short-term prediction of the system dynamics by DMD is shown in Fig. \ref{fig:pred_5415}. The observed (past) time histories are depicted in black, the predicted (future) time histories are in blue, while the true observed (future) time histories are presented with a dashed black line. All variables are shown in their standardized form and time values are normalized with the average encounter period. The NMSE of the prediction is shown at the bottom in red. Roll, pitch, and rudder present the most accurate prediction, whereas sway is found the most difficult variables to forecast. On average, variables are reasonably predicted upto two encounter periods. After, the prediction becomes less accurate especially for sway and yaw.

\subsection{KCS Turning Circle in Regular Waves}
\begin{figure}[!t]
\centering
\includegraphics[width=1\textwidth]{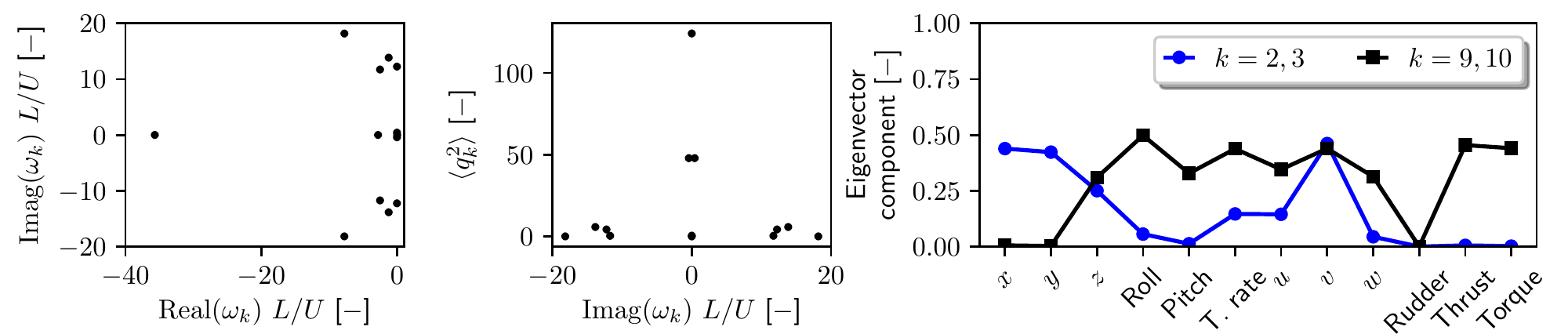}
\caption{DMD complex modal frequencies, modal participation, and two most energetic modes (from left to right, respectively) for the KCS test case.}\label{fig:freq_KCS}
\end{figure}

\begin{figure}[!t]
\centering
\includegraphics[width=0.975\textwidth]{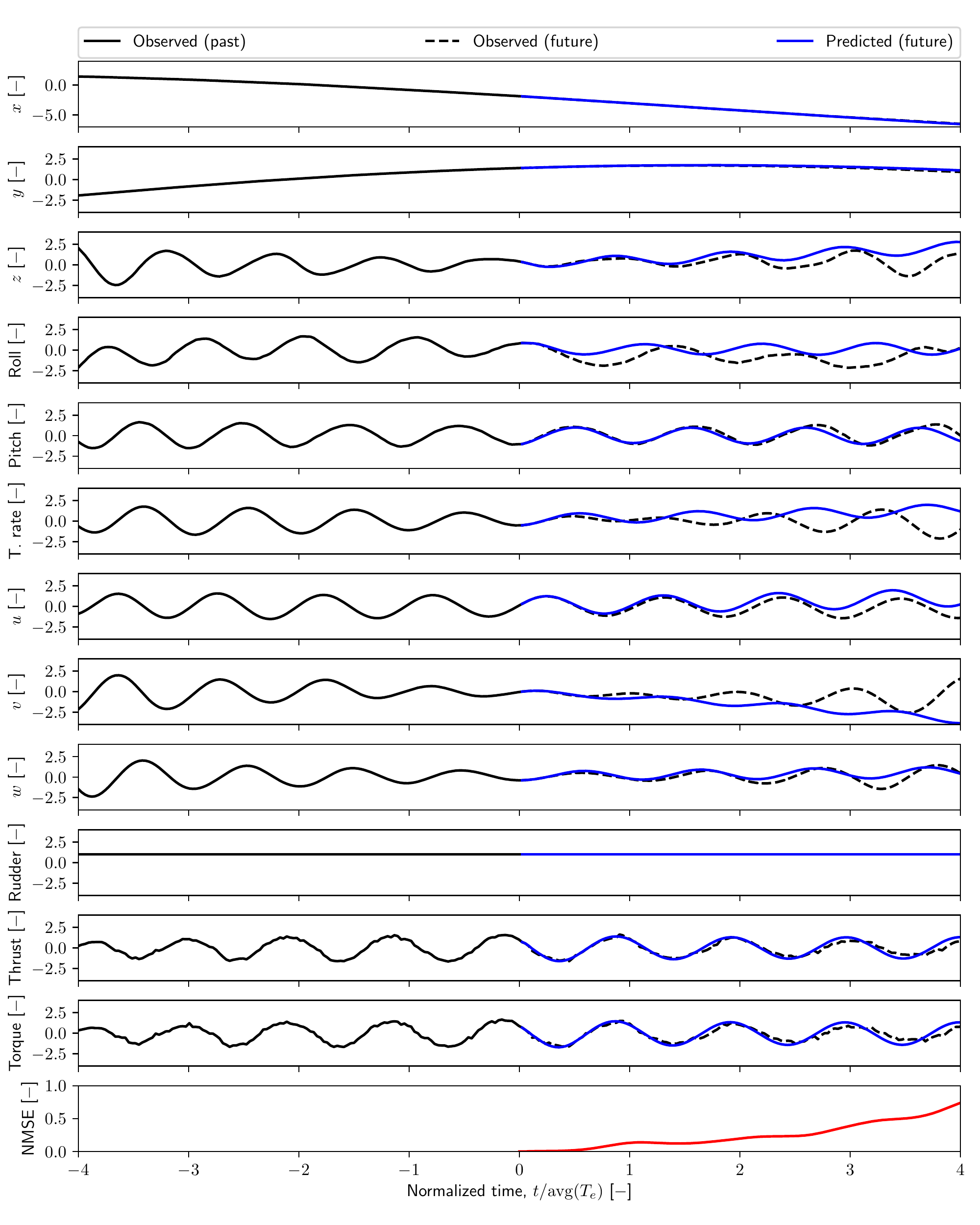}
\caption{DMD prediction and average error for KCS test case (standardized variables).}\label{fig:pred_KCS}
\end{figure}

Figure \ref{fig:freq_KCS} presents the results of the DMD of the system dynamics. Again, the left figure shows the complex modal frequencies provided by DMD. The center figure shows the modal participation as a function of the frequency (imaginary part), while the right figure presents the two most energetic couples of complex-conjugate modes (magnitude of components shown; derivatives not shown). It may be noted how the dynamic is largely participated by one mode with zero frequency (see Fig. \ref{fig:freq_KCS} center; the real part, not shown, is also zero), setting an offset for the data. The first most energetic couple of complex-conjugate modes ($k=2,3$) correspond to a quite slow dynamics and mainly involves trajectory and planar motion variables, i.e. $x$, $y$, $v$ and to a lesser extent $z$, rate of turn, and $u$. The second couple presents a larger frequency (faster dynamics, due to the waves). It mainly involves $z$, roll, pitch, rate of turn, $u$, $v$, $w$, thrust, and torque. It may be noted how the rudder does not participate in these dynamics, which may appear unexpected: in this demonstration, the rudder is kept fixed at 35\textit{deg}. Certainly, the ship dynamics depends on the rudder angle, but here the rudder is fixed (see the rudder time history in Fig. \ref{fig:pred_KCS}) and therefore no rudder dynamics is observed.

Figure \ref{fig:pred_KCS} shows the results of the DMD in producing the short-term prediction of the system dynamics, where the observed (past) time histories are depicted in black, the predicted (future) time histories in blue, the true observed (future) time histories are presented with a dashed black line. All variables are standardized and time values are normalized with the average encounter period. The NMSE of the prediction is shown at the bottom in red. The trajectory ($x$ and $y$) is very well predicted, also due to its slow dynamics. Pitch, $u$, thrust, and torque presents a faster dynamics (due to the waves) and are also reasonably predicted. Roll, turning rate, and $v$ are found the most difficult variables to forecast. Also in this case, variables are reasonably predicted upto two encounter periods. After, the prediction becomes overall less accurate.

\section{Conclusions and Future Work}
A data-driven and equation-free modeling approach for ship maneuvers in waves was presented and discussed, based on the dynamic mode decomposition of ship trajectory/motion/force time histories. The DMD provides a data-driven and equation-free approach to the linear representation of the system dynamics, allowing for (a) extracting knowledge on the system dynamics and and (b) forecasting the system's state in the near future. Results were shown for course keeping data of the self-propelled 5415M in irregular waves and turning-circle data of the self-propelled KCS in regular waves. Time histories were provided by CFD and EFD for 5415M and KCS, respectively.

Results are overall promising. The analysis is very efficient and suitable for real-time predictions. The DMD model is able to extract the most important modes and forecast with reasonable accuracy the system's state upto two wave encounter periods. After this time horizon, the prediction is no longer accurate and the methodology needs improvements.

Future research includes a systematic study of observed data size, time step, and time-derivative order along with the most effective choice of the variable set and possibly coordinate systems, including whenever appropriate dimensional analysis \cite{timur2020}. Methodological advancements that are expected to provide benefits include DMD with control \cite{proctor2016dynamic}, DMD with time delay embedding \cite{kamb2020time}, and multi-resolution DMD \cite{kutz2016multiresolution}. Finally, the combination of DMD with artificial neural network approaches \cite{dagostino2021-MARINE} is expected to overcome some of the limitations of the DMD (i.e., its linearity) providing more flexible architectures \cite{diez2022snh} to address highly-nonlinear system dynamics.

%--------------------------------------------------------------------------------
\section*{Acknowledgments}
%--------------------------------------------------------------------------------

CNR-INM is grateful to Drs. Elena McCarthy and Woei-Min Lin of the Office of Naval Research for their support through the Naval International Cooperative Opportunities in Science and Technology Program. Dr. Andrea Serani is also grateful to the National Research Council of Italy, for its support through the Short-Term Mobility Program 2018. The research is conducted in collaboration with NATO STO Research Task Groups AVT-280 ``Evaluation of Prediction Methods for Ship Performance in Heavy Weather,''  AVT-348 ``Assessment of Experiments and Prediction Methods for Naval Ships Maneuvering in Waves,'' and AVT-351 ``Enhanced Computational Performance and Stability \& Control Prediction for NATO Military Vehicles''.

\bibliographystyle{unsrt}  
\bibliography{biblio}  %%% Remove comment to use the external .bib file (using bibtex).

\end{document}